\documentclass[letter,11pt]{amsart}

\usepackage[latin1]{inputenc}
\usepackage[T1]{fontenc}
\usepackage{xy}
\usepackage{fullpage}
\usepackage{amsmath, amsthm, amssymb}
\usepackage{latexsym}
\usepackage{graphicx}
\usepackage{epsf, subfigure, verbatim, psfrag}

\usepackage{srcltx}

\newtheorem{theorem}{Theorem}[section]
\newtheorem{proposition}{Proposition}[section]

\title{A geometric refinement of a theorem of Chekanov}
\author{François Charette}
\thanks{The author holds an FQRNT doctoral research scholarship.}
\address{Department of Mathematics and Statistics, University of Montreal, C.P. 6128 Succ. Centre-Ville Montréal, QC H3C 3J7, Canada}
\email{charette@dms.umontreal.ca}

\begin{document}             

\begin{abstract}
We prove a conjecture of Barraud and Cornea (\cite{bc2}) in the monotone setting, refining a result of Chekanov on the Hofer distance between two Hamiltonian isotopic Lagrangian submanifolds. 
\end{abstract}

\maketitle

\section{Introduction}
Given $L$ and $L'$ two monotone Hamiltonian isotopic Lagrangian submanifolds (not necessarily transversal) of a symplectic manifold $(M^{2n}, \omega)$, consider an embedded symplectic ball of radius $r$ disjoint from $L'$ with real part lying on $L$.  The main result of this article relates the radius of this ball to the Hofer distance $\nabla(L,L')$ between the two Lagrangians:
\begin{theorem}\label{thmlargeur}
Let $r$ be as above, then $\frac{\pi}{2}r^2 \leq \nabla(L,L')$.
\end{theorem}

The main ingredient we will use in the proof is the following
\begin{theorem}\label{thmprincipal}
  Let $L, \; L'$ be two monotone Hamiltonian isotopic Lagrangian submanifolds of a tame symplectic manifold $M$.  Then for every almost complex structure $J \in \mathcal{J}$ and every $x \in L \backslash L'$, there exists a non constant J-holomorphic map $u$ that is either a strip with boundary on $L$ and $L'$, a disk with boundary on $L$, or a sphere such that $x \in Im(u)$ and $\int u^* \omega \leq \nabla(L, L')$.  If $L$ and $L'$ are transversal, then for a generic choice of $J$, we have $\mu(u) \leq n$.
\end{theorem}
This solves a conjecture of Barraud and Cornea \cite{bc2}, which was stated in the more general non-monotone case.
It also shows that the Hofer distance is non-degenerate, thus recovering a result of Chekanov \cite{chek1}.  By inspecting the proof of Theorem \ref{thmprincipal2}, we can rule out the presence of a sphere for a set of second category in $L\backslash L'$.

Applying Theorem \ref{thmprincipal} to a displaceable Lagrangian, we recover the well known result that through every (generic) point of $L$ there is a pseudo-holomorphic disk of symplectic area less than the disjunction energy.  This implies yet another result of Chekanov \cite{chek}, claiming that the energy needed to hamiltonialy displace a Lagrangian submanifold is at least the symplectic area of the smallest non constant pseudo-holomorphic disk or sphere.

Our theorem also recovers a result of Barraud and Cornea (\cite{bc}, \cite{bc2}) saying that when the distance between two Lagrangians is smaller than the bubbling threshold, or if there are neither pseudo-holomorphic disks nor spheres, then through every point of $L$ there is a strip of Maslov index at most n whose symplectic area is smaller than $\nabla(L, L')$.

These direct applications of Theorem \ref{thmprincipal} can be thought of as two extreme cases, one in which there are no strips, and one in which there are no disks.  The novelty of our result lies in the intermediate case, where it is a priori not clear if a holomorphic map exists.  Theorem \ref{thmprincipal} tells us that one of them exists and has a small enough energy.

The tricky part in showing such a result is that it is not enough to look at the Lagrangian Floer homology of $L$, nor its quantum homology, as they might vanish.  However, we find strips or disks of a given energy by working directly at the chain level.  The main ingredients we will use are an action of the pearl complex on the Floer complex (see \S \ref{sectionmodule}) combined with a chain homotopy between the identity and the composition PSS$^{-1} \circ$PSS (see \S \ref{sectionpss}).

We will define all the relevant structures to prove the main theorems in the next section, then proceed with their proof and give some energy estimates needed therein.  The reader familiar with Lagrangian quantum homology should go to \S \ref{sectionpreuve} and come back to the previous sections for the relevant definitions.

\section{Acknowledgements}
This research is part of my doctoral thesis under the supervision of Octav Cornea at the University of Montreal.  I thank him for his advice and guidance over the last two years.  I would also like to thank the referee for useful suggestions regarding non transversal Lagrangians.

\section{Algebraic structures}
\subsection{Definitions and conventions}\label{def}
We will only consider connected tame symplectic manifolds $(M, \omega)$ of dimension $2n$.  The set of all $\omega$-compatible almost complex structures is denoted $\mathcal{J}$ and $g_{\omega,J} ( \cdot, \cdot):= \omega( \cdot, J \cdot)$ is the associated Riemannian metric, sometimes written as a scalar product $< \cdot , \cdot>$.

Let $\mathcal{H}(M) := \{ H: M \times [0,1] \to \mathbb{R} \; | \text{ supp } H \text{ is compact} \}$ endowed with the norm
$$||H|| = \int_0^1 \max H_t - \min H_t dt.$$

The symplectic gradient of $H \in \mathcal{H}$ is the unique one-parameter family of vector fields defined by
$$\imath_{X_H^t} \omega = -dH_t.$$
These vector fields generate the Hamiltonian flow $\Psi_t$ by the differential equation
$$\frac{d}{dt} \Psi_t = X_H^t \circ \Psi_t, \; \Psi_0 = id.$$
The set of all time-1 Hamiltonian flows is called the group of Hamiltonian diffeomorphisms (or isotopies) and is denoted by $Ham(M, \omega)$.  The energy of a Hamiltonian isotopy $\phi \in Ham(M, \omega)$ is
$$E(\phi) = \inf_{ H | \Psi_1^H = \phi} \{ ||H|| \}.$$

A submanifold $L \subset M$ is Lagrangian if  $\dim L = n = \dfrac{1}{2} \dim M$ and $\omega|_{T_*L} = 0$.  We only consider closed (i.e. connected and compact) Lagrangians.  There are two important homomorphisms associated to such a submanifold, namely the symplectic area of disks with boundary on $L$:
\begin{align*}
\omega: & \pi_2(M, L) \to \mathbb{R}\\
u & \mapsto \int_{D^2} u^*\omega 
\end{align*}

\noindent and the Maslov index:
$$\mu: \pi_2(M, L) \to \mathbb{Z}.$$
The Maslov index maps a disk to the homotopy class of a loop in $\Lambda(\mathbb{R}^{2n})$, which is the set of all Lagrangian subspaces of  $\mathbb{R}^{2n}$, by following the tangent space of $L$ along the boundary of the disk.  See \cite{mcsal} for a precise definition.  The image of this morphism is then $N_L \mathbb{Z}$, where $N_L \geq 0$ is called the minimal Maslov class.

In this article, \textbf{all Lagrangians are \textit{monotone}}, which means that $N_L \geq 2$ and that there is a constant $\tau > 0$ satisfying
$$\omega(A) = \tau \mu(A) \; \forall A \in \pi_2(M,L).$$
Moreover, we will say a Hamiltonian is \textit{non-degenerate} if $L$ and $\Psi_1(L)$ intersect transversally.  As is now well known, this is a standard requirement for the definition of Lagrangian Floer homology.

Given a fixed Lagrangian submanifold $L$, let $\mathcal{L}(L) := \{ \phi(L) \; | \; \phi \in Ham \}$ be the set of all Lagrangians that are Hamiltonian isotopic to $L$.  The \textit{Hofer distance} on this set is defined by
$$\nabla: \mathcal{L} \times \mathcal{L} \to \mathbb{R}$$
$$\nabla(L, L') = \inf_{\phi \in Ham | \phi(L) = L'} \{ E(\phi)\}.$$

\subsubsection{Simple and absolutely distinct pseudo-holomorphic disks}
In order to define various algebraic structures, we need to consider spaces of pseudo-holomorphic disks with Lagrangian boundary conditions satisfying some incidence relations, e.g. the boundary should intersect some unstable manifolds transversally.  The general procedure usually goes as follows.  First, we define a space involving disks, Floer strips and other objects and compute its virtual dimension.  Second, we show that when the virtual dimension is at most one, then this space is a manifold whose actual dimension is equal to its virtual one.  The key step will often amount to showing that the disks involved are simple and absolutely distinct, which is done by an induction on the total Maslov class and relies on a technical result of Lazarrini (see \cite{bic3} and \cite{laz}).

As this has been treated thoroughly by Biran and Cornea, we will not insist on these conditions, but the reader should bear them in mind when considering the various moduli spaces.

\subsubsection{Morse theory}
Given a Morse-Smale function $f$ and $\rho$ a Riemannian metric, we denote by $Crit(f)$ the set of all critical points of $f$, $\phi_t$ the negative gradient flow and $W^u(x)$ (resp. $W^s(x)$) the unstable (resp. stable) manifold of a critical point $x$.  The Morse index of a critical point $x$ is defined by $|x| = \dim W^u(x)$ and induces a grading on $Crit(f)$.

\subsection{Lagrangian quantum homology}
\subsubsection{The pearl complex}
We recall briefly the construction of the pearl complex, first suggested by Oh \cite{oh2} following an idea of Fukaya and further developed by Biran and Cornea, whose homology is called Lagrangian quantum homology and is denoted by $QH(L)$.  The reader is invited to read the articles \cite{bic}, \cite{bic2} and \cite{bic3} for a more detailed exposition.

We denote $\Lambda = \mathbb{Z}_2 [t, t^{-1}]$ the ring of Laurent polynomials with grading given by $\deg t = -N_L$.

Let $f: L \to \mathbb{R}$ be a Morse-Smale function and $J \in \mathcal{J}$.  The \textit{pearl complex} associated to $f$ is the graded ring
$$\mathcal{C}_*(f, \rho, J) = (\mathbb{Z}_2<Crit(f)> \otimes \Lambda)_*$$
where the index of a critical point is the Morse index.

\begin{figure}[ht]
\psfragscanon
\psfrag{u1}{$u_1$}
\psfrag{uk}{$u_k$}
\psfrag{x}{$x$}
\psfrag{y}{$y$}

  \begin{center}
     \includegraphics[width=6cm]{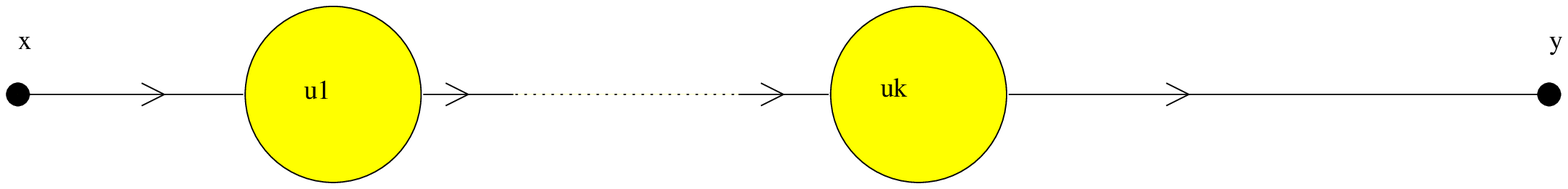}
    \caption{A pearl!}
    \label{perle}
  \end{center}
\end{figure}

The differential is defined by counting (modulo 2) Morse flow lines as well as so-called pearly trajectories, where one follows Morse flow lines connecting $J$-holomorphic disks with boundary on $L$ (see figure \ref{perle}).  This differential splits as a sum $d = \sum_k d_k$, where $d_k$ counts the number (modulo 2) of pearly trajectories with total Maslov class $kN_L$, $k \geq 0$.

To be able to count the number of pearl trajectories, one must first show that they form a manifold of dimension $0$ (with suitable index restriction).  This is shown in \cite{bic3} and relies heavily on monotonicity.  Biran and Cornea then show that $d^2 = 0$ by using a gluing argument for products of 0-dimensional pearls.  Finally, they show, by an adaptation of standard Morse cobordism arguments, that quantum homology is independent of generic choices of $J$ and $f$.

Without the monotonicity assumption, it would be necessary to enlarge the pearl complex and consider so-called clusters (see \cite{cl}), or alternatively consider the $A_\infty$ machinery of Fukaya, Oh, Ohta, Ono \cite{fooo}.

\subsubsection{The Lagrangian quantum product}
In this section, we recall how to endow quantum homology with a unitary ring structure, which coincides with the Morse-theoretic intersection product when there are no pseudo-holomorphic disks.  There are examples where this ring is non-commutative (see e.g. \cite{bic}).

\begin{figure}[ht]
\psfragscanon
\psfrag{u1}{$\mathbf{u_2}$}
\psfrag{u2}{$\mathbf{u_1}$}
\psfrag{u3}{$\mathbf{u_3}$}
\psfrag{v}{$v$}

\psfrag{x}{$y$}
\psfrag{y}{$x$}
\psfrag{z}{$z$}

  \begin{center}
     \includegraphics[width=6cm]{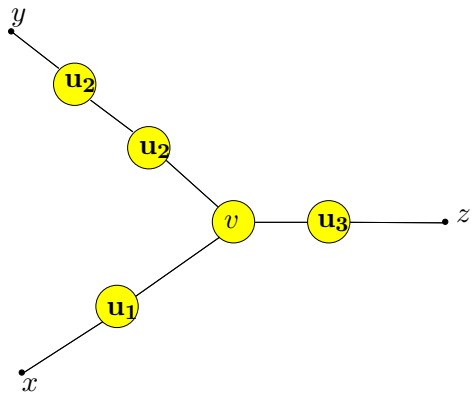}
    \caption{A tripod}
    \label{tripode}
  \end{center}
\end{figure}

Fix two Morse-Smale functions $f,g$.  Let $x \in Crit(f)$, and $y, z \in Crit(g)$.  We will consider the space of tripods of total Maslov class $kN_L$, $\mathcal{P}^{x,y}_z(k)$, as shown on figure \ref{tripode}.  The top left "pearl" leaves from $y$ and connects pseudo-holomorphic disks with flow lines of $-\nabla g$ until it reaches a possibly constant pseudo-holomorphic disk $v$.  The bottom left pearl leaves from $x$ and uses the flow of $-\nabla f$ instead, until it reaches the same disk $v$.  As for the pearl on the right, it leaves $v$ using the flow of $-\nabla g$ and eventually goes into $z$.  Notice that when no disks appear, we recover the standard Morse intersection product.  The product seen on figure \ref{tripode} is then $x \circ y = z \otimes t^k$, where $kN_L$ is the total Maslov class of the disks.

The same technical considerations as before allow us to count (modulo 2) the number of tripods of dimension 0.  This is then used to define the \textit{quantum product}
$$\circ : (\mathcal{C}(f) \otimes_\Lambda \mathcal{C}(g))_* \to \mathcal{C}(g)_{*-n}.$$

By looking at a suitable compactification of one-dimensional tripods and using a gluing argument (once again, see \cite{bic3}), one  shows that it induces a product in homology, also called the quantum product.  Standard cobordism arguments also show that the product is independent of generic choices.  

It is readily verified, by dimensional arguments, that $M_f$ is a unit on the chain level for this product, where $f$ is a Morse-Smale function having a unique maximum $M_f$, i.e. $M_f \circ y = y \; \forall y \in Crit(g)$.  However, $M_g$ would be a unit only in homology.

\subsection{Lagrangian Floer homology}
\label{subsectionlfh}
We now recall the construction of (a version of) monotone Lagrangian Floer homology (see also \cite{oh}).

Fix a non-degenerate Hamiltonian $H \in \mathcal{H}$ whose Hamiltonian flow is $\Psi_t$.  We will be interested in the set of contractible Hamiltonian orbits starting and ending on $L$.  First consider the set $\mathcal{P}_0 (L)= \{ \gamma \in C^\infty([0,1], M) \; | \; \gamma(0) \in L, \; \gamma(1) \in L, \; [\gamma] = 1 \in \pi_1(M,L) \}$ and write $\mathcal{O}_H = \{ \gamma \in  \mathcal{P}_0 (L) \; | \; \gamma'(t) = X_H^t(\gamma(t))\}$ the subset of Hamiltonian orbits therein.  As $H$ is non-degenerate, there are only a finite number of such orbits, which also correspond to the points of $L \cap \Psi_1(L)$.

Define by $D^2_- = \{ z \in \mathbb{C} \; | \; |z| \leq 1, Re(z) \leq 0\}$ the left part of the complex unit disk.  We split its boundary into two paths parametrized by $[0,1]$,  $\gamma_1 = D^2_- \cap S^1$ and $\gamma_2=D^2_- \cap i \mathbb{R}$ , oriented in such a way that $\gamma_1 \# \gamma_2$ runs  counter-clockwise.

Now let $\gamma \in \mathcal{O}_H$ and consider $u: (D^2_-, \gamma_1, \gamma_2) \to (M,L, \gamma)$ with $u(\gamma_2(t)) = \gamma(t)$, that is, $u$ is a half-disk capping $\gamma$ (as $\gamma \in \mathcal{P}_0 (L)$, such a half-disk exists).  Trivializing the symplectic bundle $u^*TM$, we obtain a path of Lagrangian subspaces of $\mathbb{R}^{2n}$ by following the path $T_{\gamma_1(t)}L \# (\Psi_t)_*(T_{\gamma_1(1)}L)$.  We may associate a Maslov index to such a path (see \cite{rs}) by first choosing $T_{\gamma_2(0)}L$ as a reference Lagrangian subspace.  With this choice, the Maslov index of a pair $(u, \gamma)$ satisfies $\mu(u,\gamma) + n/2 \in \mathbb{Z}$ (\cite{rs}, Theorem 2.4).  It is a half integer that verifies $\mu(u,\gamma) + \mu(-v, \gamma) = \mu(u \# (-v)) \in N_L \mathbb{Z}$, so it generalizes the Maslov index of a loop.

We define an equivalence relation on these pairs by  $(u,\gamma_1) \sim (v, \gamma_2) \iff \gamma_1 = \gamma_2$ and $\mu((u,\gamma_1)) = \mu((v,\gamma_2))$.  The quotient set is denoted by $\tilde{\mathcal{O}}_H$.  Notice that  $\tilde{\mathcal{O}}_H$ is in one-to-one correspondence with $\mathcal{O}_H \times \Lambda$ because there is a transitive action of $\pi_2(M,L)$ on it (given by the connected sum) with stabilizer $\ker \mu$ and  $\pi_2(M,L) / \ker \mu \cong \Lambda$.  So for each $\gamma \in \mathcal{O}_H$, we fix a representative $\tilde{\gamma}:=[u_\gamma, \gamma] \in \tilde{\mathcal{O}}_H$.

The Floer complex is defined by
$$CF(L,H,J) = \mathbb{Z}_2< \tilde{\gamma} \; | \; \gamma \in \mathcal{O}_H> \otimes \Lambda,$$
and the grading is given by $|\tilde{\gamma}| := n/2 - \mu(\tilde{\gamma})$.  With this choice of grading, the PSS and PSS$^{-1}$ morphisms (see \S\ref{sectionpss}) preserve the degree (this is not a serious issue, but it makes the notations easier to follow).

The differential is defined by counting (modulo 2, as usual) strips $u: \mathbb{R} \times [0,1] \to M$ such that $u(t, i) \in L, \; i=0,1$, satisfying \textit{Floer's equation}
\begin{align}
\label{floereq}
\partial_s u + J(u) \partial_t u  + \nabla H_t(u) = 0
\end{align}

\noindent and the asymptotic conditions
\begin{align}
\label{condas}
\lim_{s \to \pm \infty} u(s,t) = \gamma^{\pm}(t) \in \mathcal{O}_H.
\end{align}

Set
$$ \mathcal{M}(\tilde{\gamma}^-, \tilde{\gamma}^+, k)  = 
\left\{ 
u: \mathbb{R} \times [0,1] \to M \Bigg| 
\begin{array}{c}
u \; \text{verifies} \; (\ref{floereq}) \; \text{and} \; (\ref{condas})\\
\mu(u_{\gamma^-} \# u  \# -u_{\gamma^+}) = kN_L
\end{array} 
\right\}.$$

  For generic choices of hamiltonian $H$ and of $J \in \mathcal{J}$, $\mathcal{M}(\tilde{\gamma}^-, \tilde{\gamma}^+, k)$ is a manifold of dimension $|\tilde{\gamma}^-| -|\tilde{\gamma}^+| + kN_L$.  Notice that $\mathbb{R}$ acts by translation on each strip, i.e. $u(s,t) \mapsto u(s_0 + s, t), s_0 \in \mathbb{R}.$  We denote the quotient space by $\tilde{\mathcal{M}}(\tilde{\gamma}^-, \tilde{\gamma}^+, k)$.

%Finally, the Floer complex is defined by
%$$CF(L,H,J) = \mathbb{Z}_2< \tilde{\gamma} \; | \; \gamma \in \mathcal{O}_H> \otimes \Lambda,$$
%and the grading is given by $|\tilde{\gamma}| := n/2 - \mu(\tilde{\gamma})$.  With this choice of grading, the dimension of $\tilde{\mathcal{M}}(\tilde{\gamma}^-, \tilde{\gamma}^+, k)$ is $|\tilde{\gamma}^-| + kN_L - |\tilde{\gamma}^+| - 1$; moreover, as we will see in \S\ref{sectionpss}, the PSS morphism preserves the degrees 

We define the Floer differential by
$$\partial :  CF_*(L, H, J) \to CF_{*-1}(L, H, J)$$
$$\tilde{\gamma}^- \mapsto \sum_{\substack{\tilde{\gamma}^+, k  \\ |\tilde{\gamma}^-| + kN_L - |\tilde{\gamma}^+| - 1 = 0  }} \#_2 \tilde{\mathcal{M}}(\tilde{\gamma}^-, \tilde{\gamma}^+, k) \tilde{\gamma}^+ \otimes t^k$$
which we extend by linearity.

This is indeed a differential and the Lagrangian Floer homology of the complex is independent of generic $H$ and $J$ (see for example \cite{oh}).

\subsection{Lagrangian $PSS$ and $PSS^{-1}$ morphisms}\label{sectionpss}
We have now defined two homologies associated to a Lagrangian submanifold.  The next step is to recall the comparison map between them by defining suitable morphisms between the respective complexes.  It is not so surprising that we are actually computing the same homology and we recall in this section how to prove this, using the so called PSS and PSS$^{-1}$ morphisms.

These morphisms were first introduced in Hamiltonian Floer homology by Piunikhin, Salamon and Schwarz (\cite{pss}) and were then studied independently by Albers (\cite{albers} and \cite{aerr}) in the Lagrangian monotone case and by Kati\'c and Milinkovi\'c (\cite{km}) for the zero section of the cotangent bundle.  Barraud and Cornea (\cite{bc2}) also obtained some of their results by studying them in the non-monotone case under the bubbling threshold.  They appear in full generality in \cite{cl}.

\begin{figure}[ht]
\psfragscanon
\psfrag{g}{$\gamma$}
\psfrag{l}{$L$}
\psfrag{u}{$u$}

  \begin{center}
     \includegraphics[width=4cm]{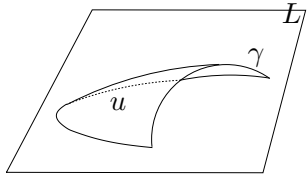}
    \caption{Half-disk for the PSS}
    \label{dpss}
  \end{center}
\end{figure}

As always, some moduli spaces will be needed.  Since we want to compare the pearl complex with the Floer complex, the geometric idea is to consider a space of pearls where the last disk is a perturbed half-disk which converges to a Hamiltonian orbit.  These perturbed half-disks can be seen as strips $u: \mathbb{R} \times [0,1] \to M$ satisfying the \textit{PSS equation}
\begin{align}
\label{psseq}
u_s + Ju_t + \beta(s)\nabla H = 0,
\end{align}
where $\beta(s)$ is 0 for $s \leq 0$,  1 for $s \geq 1$ and $\beta' \geq 0$.  Thus $u$ interpolates between the Cauchy-Riemann equation $\overline{\partial}_J u= 0$ and Floer's equation (\ref{floereq}).  We also require the following boundary and asymptotic conditions (see figure \ref{dpss}):
\begin{align}
\label{aspss}
\begin{cases}
 u(s, i) \in L, \; s \in \mathbb{R}, \; i= 0,1;\\
u(-\infty, t) = l_0 \text{ for some } l_0 \in L;\\
u(\infty, \cdot) \in \mathcal{O}_H.
\end{cases}
\end{align}

Under these conditions, $u$ can be thought of as a map $u: (D^2_-, \gamma_1, \gamma_2) \to (M, L, \gamma)$.

\begin{figure}[ht]
\psfragscanon
\psfrag{x}{$x$}
\psfrag{u1}{$u_1$}
\psfrag{u2}{$u_2$}
\psfrag{ur}{$u_r$}
\psfrag{g}{$\gamma$}

  \begin{center}
     \includegraphics[width=6cm]{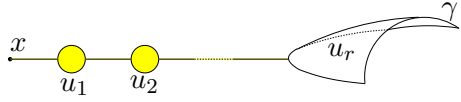}
    \caption{A pearl converging to an orbit}
    \label{pss}
  \end{center}
\end{figure}

We denote $\mathcal{M}(x, \tilde{\gamma}, m)$ the space of pearls leaving from $x$ where the last disk converges to $ \gamma$ and verifies (\ref{psseq}) and (\ref{aspss}).  The integer $mN_L$ is the total Maslov class of the disks plus the Maslov class of $-u_\gamma$ (see figure \ref{pss}).  As always, under genericity assumptions, this space is a manifold whose virtual dimension is $v:= |x| + mN_L - |\tilde{\gamma}|$ (see \cite{bic3} and \cite{albers}).  When $v \leq 1$, then $\mathcal{M}(x, \tilde{\gamma}, m)$ is a manifold of dimension $v$.
The PSS morphism is then
$$PSS: \mathcal{C}(f) \to CF(L; H)$$
$$x \mapsto \sum_{\substack{\tilde{\gamma}, m  \\|x| + mN_L - |\tilde{\gamma}| = 0}} \#_2 \mathcal{M}(x, \tilde{\gamma}, m) \tilde{\gamma} \otimes t^m.$$
It can be shown to be a chain morphism by compactifying the one dimensional spaces $\mathcal{M}(x, \tilde{\gamma}, m)$ and using a gluing argument.

The inverse morphism is defined in a similar fashion:  we start from an orbit $\gamma$ and use a half-disk to follow a pearl up to a critical point $x$.  The half-disk satisfies the \textit{$PSS^{-1}$ equation}
\begin{align} \label{pssinv}
u_s + Ju_t + \beta(-s)\nabla H = 0.
\end{align}

The space of pearls leaving from $\gamma$ using a PSS$^{-1}$ half-disk and going into $x$ is denoted $\mathcal{M}(\tilde{\gamma}, x, m)$ and its (virtual) dimension is $|\tilde{\gamma}| + mN_L - |x|$, where $mN_L$ is the total Maslov class of the disks plus the Maslov class of $u_\gamma$.  The inverse morphism is then
$$PSS^{-1}: CF(L, H) \to \mathcal{C}(f)$$
$$\tilde{\gamma} \mapsto \sum_{\substack{\tilde{\gamma}, m \\ |\tilde{\gamma}| + mN_L - |x| = 0}} \#_2 \mathcal{M}(\tilde{\gamma}, x, m)  x \otimes t^m.$$

The remaining step is to show that there is a chain homotopy $\psi : \mathcal{C}_*(f) \to \mathcal{C}_{*+1}(f)$ verifying
\begin{align} \label{pssiso}
id - PSS^{-1}PSS = (d \psi - \psi d).
\end{align}

The chain homotopy is constructed in essentially the same way as in \cite{albers}, which is adapted from \cite{pss}  in the Hamiltonian case.  We still recall how to construct it, because it will allow us to define relevant moduli spaces used in the proof of Theorem \ref{thmprincipal}.  As we will not need the fact that PSS$^{-1}$ is the right inverse of PSS, we will not prove it.

By gluing together along an orbit two half-disks verifying the PSS and PSS$^{-1}$ equations, we get an element of $\pi_2(M,L)$ satisfying a perturbed Cauchy-Riemann equation.  This disk should then be used in a one parameter family to interpolate between the identity and $PSS^{-1} \circ PSS$.  More formally, let  $u: \mathbb{R} \times [0,1] \to M$ satisfy $u(s, i) \in L, \; i=0,1,$ and the equation
\begin{align}
\label{eqalphar}
u_s + J(t, u)u_t - \alpha_R(s)\nabla H = 0, \; E(u) := \int \omega(u_s, J u_s) < \infty,
\end{align}

where $R \geq 1$ and $\alpha_R : \mathbb{R} \to \mathbb{R}$ is smooth and is such that
$
\alpha_R(s)
\begin{cases}
= 1 \text{ if } |s| \leq R \\
= 0 \text{ if } |s| \geq R+1 \\
-1 \leq \alpha_R'(s) \leq 1.
\end{cases}
$

We also set $\alpha_R = R \alpha_1$ when $R \in (0,1)$.  As $E(u) < \infty$, we can show that $u$ induces a map $u: (D^2, S^1) \to (M,L)$.  In contrast to holomorphic disks, these disks might have a negative Maslov index.

Consider next $x, y \in Crit(f)$ and a pearl going from $x$ to $y$, where one of the disks satifies (\ref{eqalphar}) for some $R$ depending on the disk and the other disks are $J$-holomorphic.  Quotient by the automorphism group and denote the set of all those disks by $\mathcal{M}(x,y, m)$, where $mN_L$ is the total Maslov class of all the disks (once again, $m$ might be negative).  The virtual dimension of this set is $v:= |x| + mN_L - |y| +1$ and when $v \leq 1$, it is a smooth manifold of dimension $v$.  The parameter $R$ increases the dimension by 1, hence the term $+1$ in the formula.

The desired chain homotopy is then given by
$$\psi:  \mathcal{C}_*(f) \to \mathcal{C}_{* +1}(f)$$
$$x \mapsto \sum_{\substack{y, m \\ |x| + mN_L - |y| +1 = 0}} \#_2 \mathcal{M}(x,y, m) y \otimes t^m.$$

To show that this is indeed a chain homotopy, we compactify spaces $\mathcal{M}(x,y, m)$ of dimension 1 and look at their boundary components.  This gives the desired property.

%The different ways a one-parameter sequence might converge to a boundary point are the following:\\
%i.  A flow line breaks on a critical point, giving rise to the terms $d \circ \psi - \psi \circ d$;\\
%ii.  The parameter $R$ diverges to infinity and breaks on two half-disks glued together along a Hamiltonian orbit.  This accounts for the term $PSS^{-1} \circ PSS$;\\
%iii. %!!!!!!!!!!ATTENTION!!!!!!!!%
%The parameter $R$ converges to $0$, thus converging to a $J$-holomorphic disk, i.e. we obtain a pearl going from $x$ to $y$, whose dimension is $0 = |x| + kN_L - |y|$.  The only possibility is $k=0$ and $x=y$, so we obtain the identity.
%ATTENTION%

\subsection{Structure of $QH(L)$-module} \label{sectionmodule}
In this section, we provide $HF(L)$ with a $QH(L)$-module structure.  This structure has been used to define spectral invariants by Schwarz (see \cite{sc}) in symplectic Floer homology (when M is symplectically aspherical) and by Leclercq (\cite{leclercq}) in Lagrangian Floer homology when there is no bubbling.

\begin{figure}[ht]
  \begin{center}
     \includegraphics[width=3cm]{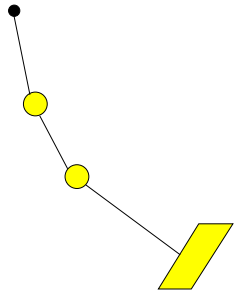}
    \caption{$x \star \gamma$}
    \label{xstargamma}
  \end{center}
\end{figure}

The geometric picture is simply to connect a pearl to a strip as on figure \ref{xstargamma}.  Let $x \in Crit(f)$ and  $\tilde{\gamma} \in \tilde{\mathcal{O}}_H$.  We then denote by $\mathcal{M}(x, \tilde{\gamma}, \tilde{\alpha}, k)$ the set of pearls leaving from $x$ such that the last flow line enters a Floer strip $v$ connecting $\gamma$ to $\alpha$ at the point $v(0,0)$, where $k$ is the total Maslov index of the pseudo-holomorphic disks, plus the Maslov class of $u_\gamma \# v \# -u_\alpha$.
 
Under genericity assumptions and when the virtual dimension is at most one, this space is a smooth manifold of dimension $|x| + |\tilde{\gamma}| + k - |\tilde{\alpha}| - n$.  We outline here a proof of this fact, as it has not been treated by Biran and Cornea.  We will only point to the references, as the details appear in various locations in the literature.  

The key idea concerns two transversality arguments.  The first one shows that unstable manifolds intersect Floer strips transversally at the point $(0,0)$.  To do this, one must allow the complex structure to vary and consider universal moduli space of strips (see \cite{fhs}).  This was used already by Leclercq \cite{leclercq}.

The second argument shows that the evaluation of a pseudo-holomorphic disk at the point $1$ intersects Floer strips at the point $(0,0)$ transversally.  Here again, one must allow the complex structure to vary and adapt the proof of \cite{mcsal2}, Proposition 3.4.2.   One then chooses a complex structure that agrees on a neighbourhood of the intersection point.  This is then used to glue together a disk and a strip into a strip on one side, and into a flow line connecting a disk and a strip on the other side.

Once these two technical steps have been taken care of, an induction argument as the one used in Proposition 3.1.3 of \cite{bic3}, shows that when the virtual dimension is at most one, all disks involved are actually simple and absolutely distinct, thus $\mathcal{M}(x, \tilde{\gamma}, \tilde{\alpha}, k)$ is a manifold of the right dimension.

The action of the pearl complex on the Floer complex is then given by
$$\star : (\mathcal{C}(f) \otimes_\Lambda CF(L, H))_* \to CF(L, H)_{* - n}$$
$$x \otimes \tilde{\gamma} \mapsto \sum_{ \substack{\tilde{\alpha}, k \\ |x| + |\tilde{\gamma}| + k - |\tilde{\alpha}| - n = 0} } \#_2 \mathcal{M}(x, \tilde{\gamma}, \tilde{\alpha}, k) \tilde{\alpha} \otimes t^k.$$

It is a chain map by the gluing argument above, and one can show it makes $HF(L)$ into a (left) $QH(L)$-module.  The proof that $(x \circ y) \star \gamma = x \star (y \star \gamma)$ in homology is a bit tricky to verify.  However, it is quite similar to the proof of (\ref{pssmod}) given below.

It has already been shown that Lagrangian quantum homology is isomorphic to Lagrangian Floer homology via the PSS morphism.  Moreover, quantum homology has an obvious $QH(L)$-module structure given by the quantum product.  It is thus natural to ask if the isomorphism preserves the module structures.  It turns out to be the case and it will be shown in the next section.

\subsection{The $PSS$ is a $QH(L)$-module isomorphism}
This property of the PSS has been shown in \cite{leclercq} when there is no bubbling.  The proof given here is a generalization.  Namely, we show that given $f, g$ two generic Morse-Smale functions, there is a chain homotopy $$\eta : (\mathcal{C}(f) \otimes \mathcal{C}(g))_* \to CF(L)_{* - n + 1}$$
satisfying
\begin{align} \label{pssmod}
PSS(x \circ y) = x \star PSS(y) + (\partial \eta - \eta \partial)(x \otimes y).
\end{align}

\begin{figure}[ht]
\psfragscanon
\psfrag{gluing}{gluing}
  \begin{center}
     \includegraphics[width=8cm]{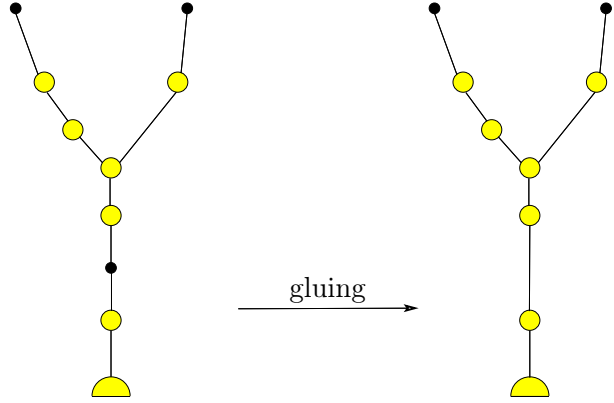}
    \caption{Gluing $PSS(x \circ y)$ to remove $z$}
    \label{pssxy}
  \end{center}
\end{figure}

We first notice that $PSS(x \circ y)$ counts the number of points (mod 2) in $\bigcup_z \mathcal{P}^{x,y}_z(k_1) \times \mathcal{M}(z, \tilde{\gamma}, k_2)$.  We can glue such a configuration along the two flow lines going into and out of $z$, thus obtaining a one parameter family of spaces that can break on the term $PSS(x \circ y)$.  Let us denote this new space by $\mathcal{S}$ for now (see figure \ref{pssxy}).

As for $x \star PSS(y)$, we can glue along the orbit $\gamma$ thus removing the Floer strip and keeping only a perturbed half-disk satisfying (\ref{pssinv}).  We denote for now by $\mathcal{S}'$ the space of such configurations.

\begin{figure}[ht]
\psfragscanon
\psfrag{glueing}{gluing}
\psfrag{bubbling}{bubbling}
  \begin{center}
     \includegraphics[width=8cm]{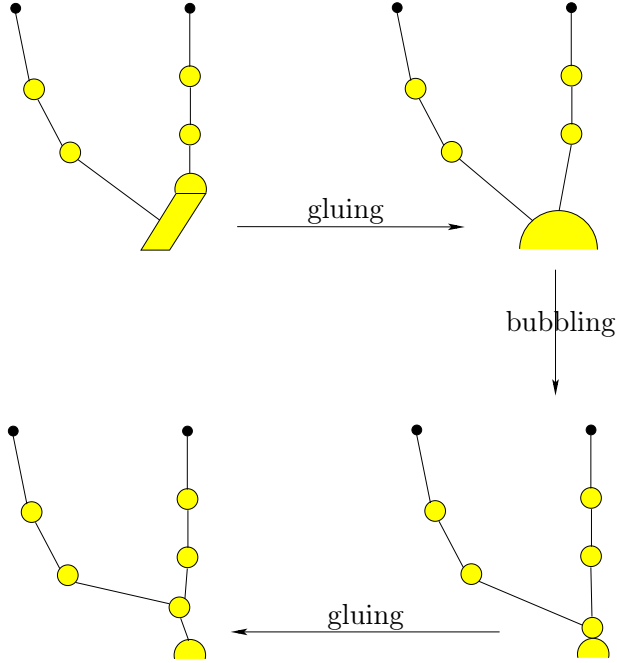}
    \caption{Gluing $x \star PSS(y)$ into $\mathcal{S}'$ and going into $\mathcal{S}$}
    \label{xpssy}
  \end{center}
\end{figure}

The point is now to use these two new spaces to define the desired chain homotopy.  Before doing so, we observe that they are actually part of the same space.  Indeed, consider a one-parameter family of elements of $\mathcal{S}'$ where a disk with two marked points bubbles out of $\mathcal{M}(H, \tilde{\alpha}, k)$ in the limit.  Using standard gluing arguments, one can then glue back this disk in the space $\mathcal{S}$. Hence the two spaces should actually be used simultaneously to define the chain homotopy.  See figure \ref{xpssy}.

This leads us to define the chain homotopy as
$$\eta: (\mathcal{C}(f) \otimes \mathcal{C}(g))_* \to CF(L)_{* - n + 1}$$
$$x \otimes y \mapsto \sum_{\substack{\tilde{\alpha}, k \\ |x| + |y| + k - |\tilde{\alpha}| - n = 0} } \#_2 \mathcal{P}^{x,y}_{\tilde{\alpha}}(k)  \tilde{\alpha} \otimes t^k$$

The space $\mathcal{P}^{x,y}_{\tilde{\alpha}}(k)$ is defined as in the quantum product, except for the last (half) disk $u$ which satisfies the PSS equation (\ref{psseq}) and converges to $\alpha$. Moreover, $kN_L$ is the Maslov class of all the disks plus the one of $-u_{\alpha}$.  Note that we allow for a pearl leaving from $x$ to end anywhere on $u((-\infty, \infty), 0)$, as in the top-right part of figure \ref{xpssy}.

Looking at spaces $\mathcal{P}^{x,y}_{\tilde{\alpha}}(k)$ of dimension 1 and their boundary combined with the previous discussion proves that $\eta$ is indeed a chain homotopy between $PSS(x \circ y)$ and $x \star PSS(y)$.

\section{Proof of theorems \ref{thmprincipal} and \ref{thmlargeur}}\label{sectionpreuve}
We have now everything needed to prove Theorem \ref{thmprincipal}.  It follows directly from the next
\begin{theorem}\label{thmprincipal2}
Let $L \subset (M^{2n}, \omega)$ be a monotone Lagrangian submanifold and $H \in \mathcal{H}$ a non-constant Hamiltonian. Then for every $J \in \mathcal{J}, \; x_0 \in L$, there exists a non constant $J$-holomorphic map $u$ that is either a Floer strip with boundary on $L$, a disk with boundary on $L$ or a sphere, such that $u$ goes through $x_0$ and $E(u) \leq ||H||$.  If $L$ and $\Psi_1^H(L)$ intersect transversely, then for a generic $J$, $0 \leq \mu(u) \leq n$.
\end{theorem}

To get Theorem \ref{thmprincipal} from this, simply apply the flow $\Psi_t^H$ associated to $H$ to the map $u$ in case it is a Floer strip, the action being given by $\Psi_t^H(u(s,t))$, thus obtaining a $J$-holomorphic strip with boundary on $L$ and the desired energy bound.

\noindent \textbf{Proof}
\textit{First case:  $\Psi^H_1(L)$ intersects $L$ transversely.}\\
Let $f,g$ be two generic Morse-Smale functions, each having a unique minimum (resp. maximum) $m_f$ and $m_g$ (resp. $M_f$ and $M_g$).  It might not be possible to choose $m_f = x_0$, the point through which we would like a strip or a disk to go.  The reason is that $W^u(x_0)$ might not intersect the space of Floer strips transversally.  However, by genericity arguments, we may take $m_f$ as close as we want to $x_0$ and then use a sequence of strips/disks converging to $x_0$ to get the result.  Notice that a sphere might bubble-off at $x_0$.  This also explains why we may not choose the Maslov index to be exactly $n$, but at most $n$ (we use here that $J$ is generic).  So we assume $m_f = x_0$.

Combining formulas (\ref{pssiso}) and (\ref{pssmod}), we get
$$ x \circ y = PSS^{-1}(x \star PSS(y)) + PSS^{-1}(\partial \eta - \eta \partial)(x \otimes y) + (d \psi - \psi d)(x \circ y)$$

Taking $x = m_f$ and $y = M_g$ yields
\begin{align*}
m_g \otimes t^0 + h.o.t & =  PSS^{-1}(m_f \star PSS(M_g))\\
& + PSS^{-1}(\partial \eta - \eta \partial)(m_f \otimes M_g)\\ 
& + (d \psi - \psi d)(m_f \circ M_g)\\
& :=(1) + (2) + (3). 
\end{align*}

The expression h.o.t stands for higher order terms and represents terms whose projection in the ring of Laurent polynomials gives polynomials of degree at least one.  Note that there might not be any such terms, but that $m_g \otimes t^0$ certainly appears on the left hand side, by standard Morse homology arguments.  Thus it also appears on the right hand side, hence in one of the expresssions (1), (2) or (3).  We conclude that one of the moduli spaces used to define these expressions must be non-empty.  We now show that in any case, we get the desired bound on the energy:

Case (1):  As $m_f$ is the minimum, $W^u(m_f)$ contains only the point $m_f$, so $m_f$ touches a Floer strip starting on $PSS(M_g)$ and ending on an orbit connected to $m_g \otimes t^0$ via PSS$^{-1}$.  Denote $u_1$ the half-disk satisfying the PSS equation, $u_2$ the Floer strip and $u_3$ the half-disk used in the PSS$^{-1}$.  Notice also that no holomorphic disks appear, for dimensional reasons.  Using the energy estimates (\ref{enerpss}), (\ref{enerfloer}) and (\ref{enerpssinv}) of \S \ref{energy}, we get
\begin{align*}
0 & < \sum_1^3 E(u_i)\\
& \leq \sum \omega(u_i) + ||H||\\
& = ||H||,
\end{align*}
where the last equality comes from the fact that the total Maslov class (hence the total symplectic area, by monotonicity) is zero, as we are considering the term $m_g \otimes t^0$.

Case (2):  First, note that $m_g \otimes t^0$ cannot be a boundary, because it is the unique minimum, so we can ignore the $\partial \eta$ part of the term (2) and consider only the $\eta \partial(m_f \otimes M_g)$ part.  Moreover, $M_g$ is a cycle by Morse theory arguments as well as for degree reasons, hence we simplify again and the only non zero term is $\eta (d m_f) \otimes M_g$.  As $m_f$ is the minimum and $d m_f \neq 0$, we conclude that a holomorphic disk goes through $m_f$.  The same energy estimates arguments as in case (1) give the energy bound.

Case (3):  The proof uses the energy estimate (\ref{enerhom}) and is identical to case (2).

\noindent \textit{Second case:  $\Psi^H_1(L)$ does not intersect $L$ transversely.}\\
This case is taken care of by considering a sequence of hamiltonian isotopies $\Psi^{H^k}_1$ converging to $\Psi^H_1$ (say of Hofer norm $\in [||H||, ||H|| + \epsilon]$) such that the intersection is transverse.  We then use the first case to get a sequence of $J$-holomorphic maps and we need to study the limiting map.  The only problematic case is that of a sequence of strips $\{u_k\}$ with right boundary on $\Psi^{H^k}_1(L)$.  As it turns out, the exact same case has already been studied by Hofer (see \cite{hofer}, Proposition 2), where it is shown that the failure of compactness is due to bubbling-off of spheres or disks.  

Hence we still get the energy bound, although we do not know how to handle the Maslov index of the limiting map, as there is no control on the ``dimension'' of $L \cap \Psi^H_1(L)$ (recall from \S \ref{subsectionlfh} the definition of the Maslov index of an orbit).
\qed

\noindent \textbf{Proof of theorem \ref{thmlargeur}}
This argument is standard and goes back to Gromov.  Consider a map $u$ as in Theorem \ref{thmprincipal}, assuming first that $u$ is either a strip or a disk.  Now let $e: (B(r), \omega_0, J_0) \to (M \backslash L', \omega, e^* J_0)$ be a symplectic embedding whose real part lies on $L$.  Taking the pull-back of $u$ by $e$ gives a $J_0$-holomorphic curve whose area, by the Schwarz reflexion principle and by the theory of minimal surfaces, is at least $\frac{\pi}{2}r^2$. By the choice of $u$, this area is also bounded above by $\nabla(L,L')$.
 
In case $u$ is a sphere, the same argument shows that $\pi r^2 \leq\nabla(L,L')$.

\qed

\begin{section}{Energy estimates} \label{energy}
As these are by now quite standard and can be found in the appendix of \cite{albers}, we provide only the statement of the energy estimates used in the proof of Theorem \ref{thmprincipal2}.  We refer to the previous sections for the relevant definitions.  The only one we recall is the definition of the energy of a continuous function $u$ depending on two variables $(s,t)$:
$$E(u) := \int \omega(u_s(s,t), J u_s(s,t)) dsdt.$$

\begin{proposition}
Given a Floer strip $u \in \mathcal{M}(\tilde{\gamma}^-, \tilde{\gamma}^+, k)$, we have
\begin{align} \label{enerfloer}
E(u) & = \omega(u) - \int H_t(\gamma^+(t))dt + \int H_t(\gamma^-(t))dt.
\end{align} 
\end{proposition}

\begin{proposition}
Let $u$ be a half-disk satisfying equation (\ref{pss}) and converging to an orbit $\gamma$, then
\begin{align} \label{enerpss}
E(u) & \leq \omega(u) - \int H_t(\gamma(t)) + \int \sup_M H_t(\cdot) dt.
\end{align}
\end{proposition}

In a similar fashion, we obtain a bound on the energy of a half $u$ disk satisfying equation (\ref{pssinv}):
\begin{align} \label{enerpssinv}
E(u) & \leq \omega(u) + \int H_t(\gamma(t)) - \int \inf_M H_t(\cdot) dt.
\end{align}

\begin{proposition}
Let $u$ be a disk satisfying equation (\ref{eqalphar}), then
\begin{align} \label{enerhom}
E(u) \leq \omega(u) + ||H||.
\end{align}
\end{proposition}
\end{section}

\end{document}